\documentclass[11pt,reqno]{amsart}
\usepackage{amssymb,latexsym}

\title{EIGENVALUES\\ OF THE TACHIBANA OPERATOR\\ 
WHICH ACTS ON DIFFERENTIAL FORMS}

\author{Sergey E. Stepanov and Josef Mike\v s}

\date{}
\begin{document}
\maketitle

\begin{center}\small\footnotesize
Department of Mathematics, Finance University 
under the Government of Russian Federation,
Leningradsky Prospect 49-55,\\  125468 Moscow, Russian Federation, 
e-mail: s.e.stepanov@mail.ru\\[1mm]

Department of Algebra and Geometry, Palacky University, 
17. listopadu 12, 77146 Olomouc, Czech Republic, 
e-mail: josef.mikes@upol.cz \\[3mm]
\end{center}

\begin{abstract}
In the present paper we show spectral properties of a little-known natural Riemannian second-order differential operator acting on differential forms.\\[3mm]

\end{abstract}

% \linenumbers

%% TeX Macros
\newtheorem{lemma}{Lemma}
\newtheorem{theorem}{Theorem}
\newtheorem{corollary}{Corollary}
\def\noi{\noindent}
\def\mg{\hbox{$(M,g)$}}
\def\d{\delta}
\def\la{\lambda}
\def\lr{\hbox{$\lambda^r$}}
\def\O{\Omega}
\def\orm{\hbox{$\O^r(M)$}}
\def\om{\omega}
\def\ph{\varphi}
\def\th{\theta}
\def\R{{\mathbb R}}
\def\H{{\mathbb H}}
\def\T{{\mathbb T}}
\def\Rc{R{\hbox{\kern-0.4em\raise0.67ex\hbox{\smash{$^\circ$}}}}\,}

\section{Introduction}
\label{sec1}

\noi{\bf 1.1.} Spectral geometry is the field of mathematics which deals with relationships between geometric structures of an $n$-dimensional Riemannian manifold  \mg\  and the spectra of canonical differential operators. The field of spectral geometry is a vibrant and active one. 

In the present paper we show spectral properties of a little-known natural Riemannian second-order differential operator acting on differential forms. 
 
The paper is organized as follows. First, we consider the basis of the space of natural Riemannian (with respect to isometric diffeomorphisms of Riemannian manifolds) first-order differential operators on differential $r$-forms $(1\leq r\leq n-1)$  with values in the space of homogeneous tensors on  \mg. This basis consists of three operators  $\{d,d^*,D\}$ where $d$ is the exterior differential, $d^*$ is the formal adjoint to $d$ exterior codifferential and $D$ is a conformal differential (see \cite{[16]}). 

Second, using basis operators, we construct the well-known the \textit{Hodge-de Rham Laplacian} 
$\Delta:=d^*d+d\,d^*$  and the little-known \textit{Tachibana operator} 
$\square:=r\,(r + 1)\, D^*D$. 

We recall that the Hodge-de Rham operator   is natural (that is, they commute with isometric diffeomorphisms of Riemannian manifolds), elliptic, self-adjoint second order differential operator acting on differential forms on  \mg.  If  \mg\  is compact, spectrums of such operator is an infinite divergent sequence of real numbers, each eigenvalue being repeated according to its finite multiplicity (see  \cite[p.~334-343]{[4]};  \cite[p.~273-321]{[5]}). Moreover, for a long time many geometers have been discussed the first eigenvalue of the Hodge-de Rham operator on a compact  \mg\ . They have different estimates its lover bounds on  \mg\  of positive curvature operator (see, for example, \cite{[4], [5], [7], [9], [12], [19], [21], [22]} and etc.).

Third, we note that the Tachibana operator $\square$ is a natural (that is, they commute with isometric diffeomorphisms of Riemannian manifolds), elliptic, self-adjoint second order differential operator acting on smooth differential forms on  \mg\  too (see \cite{[4],[18]}). In addition, we will prove that $\square$ has positive eigenvalues, the eigenspaces of $\square$ are finite dimensional and eigenforms corresponding to distinct eigenvalues are orthogonal. Moreover, we will show that the first eigenvalue of the Tachibana operator on  \mg\  with the positive curvature operator has its lover bound.
\\

\noi{\bf 1.2.} The present paper is based on our report at the International conference \textit{Differential Geometry and its Applications} (Czech Republic, August 19--23, 2013). At the same time the paper is a continuation of the research, which we started in \cite{[18]} and \cite{[19]}. 

%2. 
\section{Basic definitions and results}

\noi{\bf 2.1.}  \mg\  be an $n$-dimensional compact connected Riemannian  $C^\infty$-manifold with some metric tensor $g$.  For the metric $g$ on $M$, we denote by $\nabla$  the associated Levi-Civita connection, $R$ and $Ric$ respectively the Riemannian and Ricci curvature tensors of $g$. By  \orm\ we denote the vector space of  $C^\infty$-forms of degree $r=0, 1, 2,\dots, n$ on  \mg. When the manifold  \mg\  is oriented, we denote by  $\om_g$ the canonical $n$-form, called the \textit{volume form} of  \mg.
 
The metric $g$ induces pointwise inner products in fibers of various tensor bundles over  \mg. Then using the pointwise inner product in \orm\  and the volume form $\om_g$  we can implicitly the \textit{Hodge star operator} $*$\,: \orm\ $\to \O^{n-r}(M)$  as the unique isomorphism mapping forms of degree $r$ to forms of degree $n - r$ by the formula  $\om\wedge(*\th)=g(\om,\th)\om_g$ for all forms  $\om,\th\in$ \orm\ (see  \cite[p.~33]{[1]},  \cite[p.~203]{[14]}). 

By integrating the pointwise inner product $g(\om,\th)$  for any $\om$  and $\th$  in \orm\  we get the \textit{Hodge product} (see \cite[p.~203]{[1]})
\begin{equation}\label{(2.1)} 
\langle\om,\th\rangle=\int_M g(\om,\th)\om_g =\int_M \om\wedge*\th=\int_M *\om\wedge\th.
\end{equation}  
                              
\noi{\bf 2.2.} Bourguignon proved (see \cite{[2]}) the existence of a basis of the space of natural Riemannian (with respect to isometric diffeomorphisms) first-order differential operators on \orm\  with values in the space of homogeneous tensors on  \mg\  which consists of three operators, but only two of them to he recognized. These operators are the well-known exterior differential $d$: \orm\ $\to\O^{r+1}(M)$  and the adjoint to $d$ exterior codifferential $d^*$: $\O^{r+1}(M)\to$\orm\   which defined via the formula (see  \cite[p.~204]{[14]})
\begin{equation}\label{(2.2)}
\langle d\om,\th\rangle=\langle \om,d^*\th\rangle
\end{equation}  
for any  $\om\in\O^{r-1}(M)$, $\th\in\O^{r+1}(M)$. In addition, we have proved in \cite{[15]} that the third basis operator is the \textit{conformal differential} $D$: $\orm\to\O^1(M)\otimes\orm$  which defined by the formula    $D=\nabla-\frac1{p+1}\,d-\frac1{n-p+1}\,g\wedge d^*$.  In turn, the conformal differential $D$ admits a canonical formal adjoint operator $D^*$: $\O^1(M)\otimes\orm\to\orm$   such that     $D^*=\nabla^*-\frac1{r+1}\,d^*-\frac1{n-r+1}\,g\wedge d\circ trace$ where $\nabla^*$ is the formal adjoint operator to  $\nabla$  (see \cite{[16]}). 

From the basis  $\{d,d^*,D\}$ we can define the well-known (see  \cite[p.~34]{[1]};  \cite[p.~204]{[14]}) \textit{Hodge-de Rham Laplacian} $\Delta:=d^*d+d\,d^*$  and \textit{Tachibana strong Laplacian} (see \cite{[4], [9]})
\begin{equation}\label{(2.3)}
D^*D:=\frac1{r(r+1)}\,\left(\bar\Delta-\frac1{r+1}\,d^*d-\frac1{n-r+1}\,g\wedge d\,d^*\right)
\end{equation}                                 
where $\bar\Delta:=\nabla^*\nabla$  is a strongly elliptic operator which is known as the \textit{Bochner rough Laplacian} acting on forms (see \cite[p.~54]{[1]}; \cite[p.~78]{[13]}). 

It is well known, that the kernel of $\Delta$  consists of \textit{harmonic r-forms} (see  \cite[pp.~205-212]{[14]}). On the other hand, we have proved in \cite{[16]} that the kernel of $D^*D$ consists of \textit{conformal Killing r-forms} for all $r =  1, 2, \dots, n - 1$. 
\\

\noi{\bf Remark.} Conformal Killing \textit{r}-forms for $r\geq 2$ were introduced by Kashiwada (see \cite{[11]}) as a natural generalization of \textit{conformal  Killing  vector  fields},  which  are  also called infinitesimal conformal  transformations.  
\\

\noi{\bf 2.3.} We  give reviewing some basic properties of the Hodge-de Rham Laplacian $\Delta$, its eigenvalues and their forms which will be used later for comparison with our results. 

The Hodge-de Rham Laplacian $\Delta$ has the following properties (see \cite[p.~334-343]{[4]}; \cite[p.~273-321]{[5]}):
\begin{enumerate}
%1. 
\item A necessary and sufficient condition for $\Delta\om=0$  is that  $d\om=d^*\om=0$. 
%2.  . 
\item $*\Delta=\Delta*$. If $\om$ is a harmonic form, so is  $*\om$.
%3.   
\item
$\Delta$ is self-adjoint, that is, $\langle\Delta\om,\th\rangle=\langle\om,\Delta\th\rangle$  for all  $\om,\th\in\orm$.
\end{enumerate}

We denote by $\H^r(M,\R)$ the vector space of harmonic $r$-forms defined on  \mg. Then (see \cite[pp.~202-212; 375-392]{[14]}) for each integer $r$ with $0\leq   r\leq   n$, 
$\dim \H^r(M,\R)=   b_r(M)$ for the Betti number $b_r(M)$ of   \mg. The Betti numbers $b_0(M)$, $b_1(M)$, $\dots$, $b_n(M)$ are well known topological invariants of $M$. These numbers satisfy the following Poincare duality property $b_r (M)  =  b_{n - r} (M)$.

A real number  \lr\ for which there is a $r$-form  $\om$ which is not identically zero such that 
$\Delta\om=\lr\om$  is called an \textit{eigenvalue} of  $\Delta$ and the corresponding $r$-form  $\om$ an \textit{eigenform} of $\Delta$  corresponding to  \lr. The eigenforms corresponding to a fixed \lr\  form a subspace of  \orm, namely the \textit{eigenspace} of  \lr.
 
The following statements about eigenvalues of   and their corresponding forms are valid (see \cite[pp.~334-343]{[4]}; \cite[pp.~273-321]{[5]}).
\begin{enumerate} 
%1.	
\item The Laplacian $\Delta$ has a positive eigenvalue \lr\ and in fact a whole sequence of eigenvalues which diverge to $+\infty$. 
%2.	  
\item $\lr=\la^{n-r}$ for all $r=    1,\dots, n-1$.
%3.	
\item The eigenspaces of $\Delta$  are finite dimensional. 
%4.	
\item The eigenforms corresponding to distinct eigenvalues are orthogonal.
%5.	
\item Let $\bar R$: $\O^2(M)\to\O^2(M)$  be the curvature operator of  \mg, which is uniquely defined by the following identity 
$g(\bar R(X\wedge Y),W\wedge Z) = g(R(X\wedge Y),Z,W)$
 for $X,Y,Z,W\in C^\infty TM$ (see  \cite[p.~36]{[14]}). If $\bar R$ satisfies the inequality     $\bar R\geq \d$  for some positive number  $\d$ at every point of $M$, then 
  $$\lr\geq \inf\{r(n-r+1)\,\d; (n-r)(r+1)\,\d\}  $$ 
and the Betti numbers $b_r(M)=0$ for all $r=1,\dots, n - 1$ (see \cite[pp.~342-343]{[4]},  \cite{[22]}). The equality is attained on for the Euclidian sphere of constant curvature $\d$  and some conformal Killing $r$-form (see \cite{[22]}).
\end{enumerate} 

\noi{\bf 2.4.} Let  \mg\  be an $n$-dimensional compact connected Riemannian manifold. By \cite{[15]}, the form  $\om\in\orm$ is \textit{conformal Killing} if and only if $D^*D\om=0$. It is clear that conformal Killing forms of a given degree will form a linear vector space (with constant real coefficients) over  \mg. We denote by $\T^r(M, \R)$ the vector space of conformal Killing $r$-forms defined on  \mg. Then for each integer $r$ with $1\leq r\leq n-1$, $\dim \T^r(M, \R) =  t_r(M) <\infty$   for the Tachibana number $t_r(M)$ of   \mg\  and all $r = 1,\dots, n - 1$ (see \cite{[17],[18]}). The Tachibana numbers $t_1(M),\dots, t_{n - 1}(M)$ are conformal scalar invariants of  \mg\  and satisfy the following duality property $t_r (M) = t_{n - r} (M)$ that is an analog of the Poincare duality for the Betti numbers (see \cite{[17]}). 

We define the \textit{Tachibana operator} by the formula $\square=r(r+1)\,D^*D$   and formulate some basic properties of the Tachibana operator $\square$ in the following lemma.
\begin{lemma}% lemma 1
 Let  \mg\  be an n-dimensional compact connected Riemannian manifold and \orm\  be the vector space of  $C^\infty$-forms of degree $r =  1, 2,\dots, n-1$ on  \mg. Then the Tachibana operator $\square{:}$  \orm\ $\to$ \orm\  has the following properties:
 \begin{enumerate}
%1. 
\item $\square$ is an elliptic and self-adjoint operator.
%2.   
\item $*\,\square=\square\,*$. If  $\om$ is a conformal Killing form, so is  $*\om$.
%3. 
\item The kernel of $\square$ is the finite dimensional vector space $\T^r(M,\R)$ of conformal Killing r-forms.
\end{enumerate}
\end{lemma}

A real number \lr\  for which there is some $r$-form $\om$  which is not identically zero such that $\square\,\om=\lr\om$  is called an \textit{eigenvalue} of $\square$ and the corresponding $r$-form $\om$  an \textit{eigenform} of $\square$ corresponding to \lr\ (see \cite{[19]}). The eigenforms corresponding to a fixed \lr\  form a subspace of  \orm, namely the eigenspace of  \lr.

We formulate some basic properties of eigenvalues of the Tachibana operator $\square$ and their forms in the following theorem, which will be proved later.
\begin{theorem} %Theorem 1.  
Let  \mg\  be an n-dimensional compact connected Riemannian manifold and $\square$ be the Tachibana operator. Then the following statements about eigenvalues of $\square$ and their corresponding forms are satisfied. 
 \begin{enumerate}
 %1.	
 \item The Tachibana operator $\square$ has positive eigenvalues  \lr\ and $\lr=\la^{n-r}$  for all $r=1,\dots, n - 1$.
%2.	
\item The eigenspaces of $\square$ are finite dimensional. 
%3.	
\item The eigenforms corresponding to distinct eigenvalues are orthogonal.
\end{enumerate}
\end{theorem}

We denote by $S^r(M)$ the vector space of (smooth) covariant symmetric $r$-tensor fields on  \mg. The curvature tensor $R$ defines a linear endomorphism  $\Rc $: $S^2(M)\to S^2(M)$, which called as the \textit{curvature operator of second kind} (see \cite[p.~51-52]{[1]}; \cite{[3],[20]}) by the formula
$$
(\Rc \ph)(X,Y)=\sum_{i=1}^{n}\ph(R(x,e_i),Y,e_i)
$$
for any $\ph\in S^2(M)$  and  $X,Y\in C^\infty TM$, where $\{e_1,\dots,e_n\}$  is an orthonormal basis of $T_xM$ at an arbitrary point  $x\in M$.

The symmetries of $R$ imply that  $\Rc $ is a selfadjoint operator, with respect to the pointwise inner product on  $S^2(M)$. Hence the eigenvalues of $\Rc $  are all real numbers at each point  $x\in M$. We note that if the eigenvalues of $\Rc $  are  $\geq\d$ (respectively  $\leq\-\d$) at an arbitrary point $x\in M$  then the sectional curvatures at $\geq\d$  are   (respectively  $\leq\-\d$).

It is well-known that $S^2(M)$  is not irreducible at each point  $x\in M$. We denote by $S^2_0(M)$  the vector space of traceless symmetric 2-tensor fields on  \mg. Then $S^2(M)$  splits into  $O(n,\R)$-irreducible subspaces as   $S^2(M)=S^2_0(M)\otimes C^\infty M\cdot g$     (see \cite[p.~46]{[1]}). Thus, we say $\Rc $  is positive (respectively negative), or simply   $\Rc > 0$ (respectively  $\Rc < 0$), if all eigenvalue of $\Rc $  restricted to $S^2_0(M)$  are positive (respectively negative). It should be noted that $\Rc $  does not preserve the subspace $S^r_0(M)$  in general, but it does, for instance, when $g$ is an Einstein metric (see \cite[p.~52]{[1]}). 

We can formulate now the following theorem. 
%Theorem 2. 
\begin{theorem}
Let  \mg\  be an n-dimensional compact connected Riemannian manifold and $\square$ be the Tachibana operator. Suppose the curvature operator of the second kind $\Rc $\,:  $S^2(M)\to S^2(M)$ is negative and bounded above by some negative number -- $\d$   at each point  $x\in M$, then an arbitrary eigenvalue of  $\square$  satisfies the inequality  $\lr\geq r(n-r)\,\d$ and the Tachibana numbers $t_r(M)=0$ for all $r=1,\dots, n - 1$. The equality $\lr=r(n-r)\,\d$  is attained for some harmonic eigenform of  $\square$ and in this case the multiplicity of \lr\  less or is equal to the Betti number $b_r(M)$.
\end{theorem}
\noi{\bf Remark.} It is clear that the multiplicity of the eigenvalue  $\lr=\inf\{r(n-r+1)\,\d;(n-r)(r+1)\,\d\}$   of the Hodge-de Rham Laplacian $\Delta$  less or is equal to the Tachibana number $t_r(M)$ for all $r = 1, 2,\dots, n - 1$.
\\

Suppose now that $(\H^n, g_0)$ is a compact $n$-dimensional hyperbolic manifold with standard metric $g_0$ having constant sectional curvature equal to --\,1. In this case, we can obtain from above theorem the following corollary.
\begin{corollary}
 Let $(\H^n, g_0)$ be a compact n-dimensional hyperbolic manifold then  $\lr\geq   r (n - r)$.  The equality $\lr= r (n - r)$  is attained if and only if $n=2r$. In this case the multiplicity of \lr\  is equal to the Betti number $b_r (\H^{2r})$. 
\end{corollary}

%3. 
\section{Proofs of statements}

In the remainder of this paragraph, we assume that  \mg\  is compact. Taking its orientable double covering, if necessary, we may consider  \mg\  as orientable without loss of generality.
\\

\noi{\bf 3.1.} First of all we  prove the Lemma. We know that $D^*D$ is an elliptic operator, its kernel has a finite dimension and that an arbitrary form $\om$  is \textit{conformal Killing} if and only if $D^*D=0$ (see \cite{[16]}). Therefore the Tachibana operator $\square$ is elliptic too, its kernel consists of conformal Killing forms and has the finite dimension $t_r(M)$. Next, thanks to (\ref{(2.3)}), we have
\begin{equation}\label{(3.1)}         
\langle\square\om,\th\rangle= r(r+1)\,\langle D^*D\om,\th\rangle=
r(r+1)\,\langle D\om,D\th\rangle=r(r+1)\,\langle \om,D^*D\th\rangle=\langle \om,\square\th\rangle
\end{equation}        
for any   $\om,\th\in C^\infty\Lambda^rM$, it follows that that $\square$ is a self-adjoint positive semi-definite operator. Therefore all eigenvalues of $\square$  are real positive numbers. Moreover, the equality   $*\,\square=\square\,*$   holds because for an arbitrary form $\om\in\orm$  we have
$$
*\,(\square\om)=*\,(\nabla^*\nabla\om)-\frac1{r+1}\,*(d^*d\om)-\frac1{n-r+1}\,*(d\,d^*\om)=
$$
$$
=\nabla^*\nabla(*\om)-\frac1{r+1}\  d\,d^*(*\om)-\frac1{n-r+1}\ d^*d(*\om)=
\square(*\om),
$$
where  $*\,\nabla=\nabla\,*$, $*d^*d=d\,d^**$   and   $*d\,d^*=d^*d*$ (see \cite[p.~72]{[8]}; \cite[p.~107]{[10]}; \cite[p.~9]{[23]}). The Lemma is proved.
\\

Next we  prove Theorem 1. First, from the Lemma we can conclude that the Tachibana operator $\square$  has positive eigenvalues \lr\  and $\lr=\la^{n-r}$  for all $r=1,\dots, n - 1$. 

Second, for an eigenvalue \lr\  of $\square$ we can define the differential operator of second order $\square'$:  \orm\ $\to$ \orm\ such as $\square'=\square-\lr\,\rm Id$. It is obvious that the operator $\square'$  is elliptic, its kernel has a finite dimension and consists of eigenforms of $\square$ corresponding to \lr. Therefore, the eigen-space of $\square$ corresponding to \lr\ is finite dimensional. 
 
Third, let us put $\square\om=\mu^r\om$   and  $\square\th=\lr\th$     for $\om,\th\in C^\infty\Lambda^rM$  and   $\lr\neq\mu^r$. Then it follows from (\ref{(3.1)}) that  
$\mu^r\langle\om,\th\rangle=r(r+1)\langle D\om,D\th\rangle=\lr\langle\om,\th\rangle$.
 As a corollary we obtain  $\langle\om,\th\rangle=0$. The converse is evident. This completes the proof of Theorem 1.
 \\
 
\noi{\bf 3.2.} Finally, we  prove Theorem 2. Let  \mg\  be a compact Riemannian manifold. First of all, we note that from the definition of $\square$ and (\ref{(2.3)}) we obtain 
\begin{equation}\label{(3.3)}
\langle\square\om,\om\rangle=\langle\bar\Delta\om,\om\rangle-\frac1{r+1}\ \langle d\om,d\om\rangle-\frac1{n-r+1}\ \langle d^*\om,d^*\om\rangle.
\end{equation}
for any $r$-form $\om$  on  \mg. 

In turn, the Laplacian $\Delta$  admits the following well known \textit{Weitzenb\"ock decomposition} (see \cite[p.~53]{[1]}; \cite[p.~211]{[14]}): $\Delta\om=\bar\Delta\om+F_r(\om)$  where $F_r(\om)$  depends linearly on the curvature $R$ and Ricci $Ric$ tensors of  \mg\  and  
$\langle\Delta\om,\om\rangle=\langle d\om,d\om\rangle+\langle d^*\om,d^*\om\rangle\geq0$. In this case the equation (\ref{(3.3)}) can be rewritten in the form 
\begin{equation}\label{(3.4)}
\langle\square\om,\om\rangle=-\langle F_r(\om),\om\rangle-\frac{r}{r+1}\ \langle d\om,d\om\rangle-\frac{n-r}{n-r+1}\ \langle d^*\om,d^*\om\rangle.
\end{equation} 

Next, we suppose that the curvature operator of the second kind $\Rc $: $S^2(M)\to S^2(M)$  is negative and bounded above by some negative number --\,$\d$   at each point  $x\in M$, then (see \cite{[20]}) 
\begin{equation}\label{(3.5)}
g(F_r(\om),\om)\leq-r(n-r)\,\d\,g(\om,\om).
\end{equation}                                                     
If we account of (\ref{(3.4)}), we will have 
\begin{equation}\label{(3.6)}
\langle\square\om,\om\rangle\geq r(n-r)\,\d\,g(\om,\om)+\frac{r}{r+1}\ \langle d\om,d\om\rangle+\frac{n-r}{n-r+1}\ \langle d^*\om,d^*\om\rangle.
\end{equation} 
Then for an eigenform $\om$  corresponding to an eigenvalue \lr, (\ref{(3.6)}) becomes the inequality
\begin{equation}\label{(3.7)}
(\lr-r(n-r)\,\d)\,\langle\om,\om\rangle\geq \frac{r}{r+1}\ \langle d\om,d\om\rangle+\frac{n-r}{n-r+1}\ \langle d^*\om,d^*\om\rangle.
\end{equation}                      
which proves that
\begin{equation}\label{(3.8)}
\lr\geq r(n-r)\,\d>0.
\end{equation}   
If the equality is valid in (\ref{(3.8)}), then from (\ref{(3.7)}) we obtain  $d\om=d^*\om=0$. In this case $\om$  is harmonic and the multiplicity of \lr\  less or is equal to the Betti number $b_r(M)$.
 
In particular, let  \mg\  be a compact model hyperbolic space $(\H^{2r}, g_0)$ with standard metric $g_0$ having constant sectional curvature equal to --\,1 then  $\lr\geq r (n - r)$. At the same time it is well known (see \cite{[6]}) that $L^2$-harmonic $r$-forms appear on a simply connected complete hyperbolic manifold  \mg\  of constant sectional curvature --\,1 if and only if $n=2r$. Therefore, if  \mg\  is a compact model hyperbolic space $(\H^n, g_0)$ then the equality $\lr=r(n-r)$  is attained if and only if $n=2r$. In this case the multiplicity of $\lr$  is equal to the Betti number $b_r(\H^{2r})$.

%% References


%\bibliographystyle{elsarticle-num}



\end{document}